\documentstyle[leqno]{article}
\makeatletter\@addtoreset{equation}{section}\makeatother

\begin{document}
\begin{center}
{\large\bf AN ANALYTIC SOLAR MODEL}\\
{\bf PHYSICAL PRINCIPLES AND
MATHEMATICAL STRUCTURE}\\[0,5cm]
{\bf Hans J. Haubold}\\
{\it\small UN Outer Space Office, Vienna International
Centre, A-1400 Vienna,
Austria}\\
\end{center}
\par
\bigskip
\noindent
\begin{center}
{\bf Abstract}
\end{center}
{\it The physical principles and the mathematical structure involved in
deriving an analytical representation of the internal structure of the
sun are discussed. For a two-parameter family of a non-linear matter
density distribution, the run of mass, pressure, temperature, and
luminosity throughout the sun are expressed in terms of Gauss'
hypergeometric function. The system of differential equations
governing hydrostatic equilibrium and energy conservation for the
sun generates another field of application of
special functions.}

\section{Introduction}
\hspace{\parindent}The structure of the sun is determined by
conditions associated with mass
conservation, momentum conservation, energy conservation, and
specific modes of energy transport through the sun. When considering its internal structure, rotation and magnetic
fields can be neglected since the sun is spherically symmetric.
It may come as a surprise that much of what has been observed and theorized about
the sun can be accounted for in terms of very basic physical laws:
Newton's laws of gravity and motion, the first two laws of
thermodynamics, Einstein's law of the equivalence of mass and
energy, Boyle's law, Charles' law of perfect gases, and
Heisenberg's uncertainty principle. The outline of the theory of
the structure and evolution of the sun has been formulated in  the
first half of the twentieth century; it is connected with names
like Lane, Emden, Schwarzschild, Eddington, Chandrasekhar, Hoyle,
and Fowler, see Mathai and Haubold (1988). In the second half of this
century, theories about the sun were greatly refined; this is due in part to new observation
techniques and to computer simulations of its structure and
evolution. By and large, the overall picture of the structure and
evolution of the sun seem to be well understood. However, there is
a serious discrepancy between the theory of how the sun shines and
the most direct experimental test of this theory. This discrepancy
is more specifically called 'the solar neutrino problem'; it refers
to the fact that the sun is a volume source of neutrinos which are particles
produced by thermonuclear reactions in the deep interior of the
sun, and that the copious flux of solar neutrinos predicted by
theory does not match the flux detected by experiments on earth over
the past 25 years, see Abdurashitov et al. (1994), Anselmann et al. (1994), Davis (1993), and Nakamura (1993). This problem has been studied by  Davis and Bahcall
among others, see Bahcall (1989). Solar neutrinos are the only particles that have
the ability to travel from the center of the sun to its surface
almost without interaction with solar matter, escaping freely into
space carrying the most direct information about physical processes
in the deep solar interior. The solar neutrino problem, which shall
not be the subject of this paper, stimulated further studies on
solar models, both by employing modern computing tools and using
the analytical techniques of mathematics. This was also the
justification for an effort to reconsider the derivation
of analytic solutions to the system of differential equations of
solar structure based on the very basic physical laws mentioned
above. Particularly, the solar neutrino problem has been considered to be a
reason to pursue more actively the problem of obtaining analytic formulae which give a description of the gravitationally stabilized solar fusion
reactor, thus showing that methods of the integration theory of
generalized special functions applied to solar physics constitute
another field of application for these functions, see Mathai and
Haubold (1988).\par
The nuclear reactions which cause the sun to evolve are
sufficiently
slow that the sun may be assumed to pass through a series of
equilibrium configurations. The model for the current internal structure of
the sun may be thought of as representing the sun at an
instant of time. Separating the time dependence of the evolution of
the sun from the equations governing its internal structure allows one 
to replace the time dependent partial differential equations by
four simultaneous, non-linear, ordinary differential equations of
the first order. The four equations represent the radial gradients
of mass, $M(r)$, pressure, $P(r)$, temperature $T(r)$, and
luminosity, $L(r)$, where $r$ denotes the radial distance from the 
centre of the sun. Since there are four equations but more than
four unknown
physical variables, one needs additional constitutive equations,
before the system can be solved: an equation of state for solar
matter, a nuclear energy production rate, and an opacity law. The
full system of equations must be solved subject to at most four
boundary conditions at the surface and the centre of the sun. These
boundary conditions ensure that the structure of the sun can be
calculated from the four differential equations but they do not
ensure that there is a single unique solution, see Chandrasekhar (1939)
and Stein (1966). The
procedure of numerically integrating the solar structure equations
takes advantage of large electronic computers and makes it possible
to include a variety of detailed physical effects and to vary
parameters at will, see Bahcall (1989) and Noels et al. (1993). A second procedure 
providing a solar model, whose contents can be understood intuitively without resorting to
numerical techniques, is based on Buckingham's theorem which states that a
system characterized by n physical variables can be described by a
set ensemble of n-r dimensionsless products of variables, where $r$ is
the
number of variables whose dimensionless representations are
linearly independent. This approach provides a qualitative
explanation of the fundamental stellar structure equations through
dimensional analysis. This analysis suggests that more detailed
physical and mathematical theories are essentially theories of
factors of proportionality. They eventually yield numerical values
for these proportional factors because more physical assumptions
have to be made, see Bhaskar and Nigam (1991). The third procedure
treats the solar structure
equations by rigorous mathematics leading to the Lane-Emden
equation which is a second-order non-linear differential equation
describing the structure of a polytrope gas sphere. However,
explicit analytic solutions of the Lane-Emden equation exist only
for the values n=0,1, and 5 of the polytrope of index n, not covering
the specific physical model for the internal structure of the
sun, see Chandrasekhar (1939) and Horedt (1990).\par
In addition to those three procedures for constructing solar models,
there exists another approach for finding solutions to the
solar structure equations which consists of making a specific assumption for
straightforward analytical integration of these equations. It is
possible to obtain analytic solar models by separating the
hydrostatic component from the energy-transport component of the
structure equations. For that purpose, an analytic density
distribution assumption, namely, that the matter density in the sun varies
non-linearly from the center to the thought surface, where the
density goes to zero, must be made see Stein (1966) and Mathai and
Haubold (1988).
Itis then possible to integrate the equations of mass
conservation, hydrostatic equilibrium, and energy conservation
through the sun. Together with the equation of state of a perfect
gas, the run of density, $\rho(r)$, mass, $M(r)$, pressure, $P(r)$,
temperature, $T(r)$, and luminosity, $L(r)$, are determined and can
be derived in the form of analytic formulae. The physics of the
problem requires only three independent boundary conditions:
$M(r)\rightarrow 0$ and $L(r)\rightarrow 0$ at radial
distance
$r=0;\, T\rightarrow T_0=0$ and $\rho \rightarrow \rho_0=0$
at
the radius $r=R_\odot$ of the gaseous configuration. The
requirement that
$\rho$ and $T$ tend simultaneously to specific values, in this case
zero, is only one condition since the point at which this occurs is
arbitrary. This ambiquity can be removed, in principle, by
assigning the total mass. The boundary is required to be at the
point where $M(r=R_\odot) =M_\odot$ and this provides the fourth
condition.
Hence, the central density, pressure, temperature, and total rate
of energy generation are determined as a function of the sun's
mass and radius. However, by assuming an analytic matter density
distribution, the energy-transport equation of the system of
structure differential equations can be satisfied at only one
typical point in the sun. The procedure thus established to
construct an analytic model of the solar interior allows to
determine the factors of proportionality which remain to
be an open problem in the dimensional analysis. The procedure also
reveals that the run of all physical variables for the solar model
can be expressed in terms of Gauss' hypergeometric function. 

\section{Matter Density Distribution}
\hspace{\parindent}For the integration of the system of
differential equations
governing the internal structure of the sun one has to make a
choice for an unknown function that still leaves room for physical
justification of this choice. By intuition one expects that  mass
is an increasing function while density, pressure and temperature
are
decreasing functions from the centre to the surface of the sun. Thus we make the working
hypothesis that the matter density distribution $\rho(r)$ varies
with the distance variable $r$, as
\begin{equation}
\rho(r)=\rho_c\left[1-\left(\frac{r}{R_\odot}\right)^\delta
\right]^\gamma \,, \delta > 0,\,
\gamma > 0,\, 0\leq \frac{r}{R_\odot}\leq1,
\end{equation}
where $\delta$ and $\gamma$ are kept as free parameters to ensure
that the density distribution determines properly
the mass, pressure, and temperature distributions in 
the sun. Equation (2.1) takes into account that the chosen density
distribution reflects the central value of the density
$\rho(r=0)=\rho_c$ and satisfies the boundary conditon
$\rho(r=R_\odot)=0$, where $R_\odot$ denotes the solar radius.
Also, equation (2.1) implies that $\rho\propto M_\odot/R_\odot^3$, where the
constant
of proportionality depends only on the radial mass distribution and the
radial distance.

\section{Distribution of Mass}
\hspace{\parindent}If $M(r)$ represents the total mass contained
within
the radius $r$, and
$\rho(r)$ is the density at $r$, then
\begin{equation}
\frac{dM(r)}{dr}=4\pi r^2\rho (r).
\end{equation}
Using the assumed non-linear density distribution in equation
(2.1),
integration of (3.1) throughout the Sun leads to the distribution
of
mass
\begin{equation}
M\left(\frac{r}{R_\odot}\right)=\frac{4\pi}{3}\rho_c R_\odot^3
\left(\frac{r}{R_\odot}\right)^3 \, _2F_1\left(-\gamma,
\frac{3}{\delta};\frac{3}{\delta}+1;\left(\frac{r}{R_\odot}
\right)^\delta \right),
\end{equation}
where $_2F_1(.)$ denotes Gauss' hypergeometric function, which contains
the parameters $\delta$ and $\gamma$ of the matter density
distribution in equation (2.1), see Luke (1969) and Mathai (1993). Equation
(3.2) satisfies the boundary
condition $M(r=0)=0$ and can be used to determine the central value
of the matter density distribution $\rho_c$ in equation (2.1) in
terms of the parameters $\delta$ and $\gamma$ of the chosen model
of the sun. The condition $M(r=R_\odot)=M_\odot$ in equation (3.2)
reveals that
\begin{equation}
\rho_c=\frac{3M_\odot}{4\pi
R_\odot^3}\frac{1}{\gamma!}\prod^\gamma_{i=1}\left(\frac{3}{\delta}+i\right),
\end{equation}
if $\gamma$ in equation (2.1) is kept as a positive integer; (3.3) is obtained by using the following relation for Gauss' hypergeometric
function of argument one, see Luke (1969) and  Mathai (1993):
\[_2F_1(a,b;c;1)=\Gamma(c)\Gamma(c-a-b)/\Gamma(c-a)\Gamma(c-b).\]
Equation (3.3) can be used to select the appropriate values of the
parameters $\delta$ and $\gamma$ specifying the solar model with
matter density distribution given in equation (2.1).

\section{Distribution of Pressure}
\hspace{\parindent}If $g=GM(r)/r^2$ is the gravitational force per
unit mass at $r$
due
to the attraction of the mass interior to $r$, then
\begin{equation}
\frac{dP(r)}{dr}=-\frac{GM(r)\rho(r)}{r^2}
\end{equation}
is the equation of hydrostatic equilibrium of the spherical self-
gravitating sun with $dP(r)/dr$ being
the pressure gradient. The internal pressure produced by the weight
of the overlying layers increases towards the centre while the gas
and radiation pressure must increase correspondingly to achieve the balance of
forces for equilibrium. This increase is obtained by inward
increases of both temperature and density. The internal pressure
needed to achieve a balance is the gravitational force per unit mass
$(GM_\odot/R_\odot^2)$ times the mass per unit area $(M_\odot/R_\odot^2)$ 
which gives that
$P\propto GM_\odot^2/R_\odot^4$ for any spherical body in hydrostatic self-
gravitation, where the constant of proportionality is again
determined by the radial distribution of mass in the sun and the
particular radial distance at which $P$ is measured. The constant
of proportionality can be determined by integrating equation (4.1)
throughout the volume of the sun by using equations (2.1) and (3.2) for the
density
and mass distribution, respectively. We obtain
\begin{eqnarray}
P(\frac{r}{R_\odot})& = & 
\frac{9}{4\pi}G\frac{M^2_\odot}{R^4_\odot}\left[\frac{1}{\gamma!}\prod^\gamma_{i=1}\left(
\frac{3}{\delta}+i\right)\right]^2
\nonumber \\
& & \times \frac{1}{\delta^2}\sum
^\infty_{m=0}\frac{(-
\gamma)_m}{m!(\frac{3}{\delta}+m)(\frac{2}{\delta}+m)}\left[\frac
{\gamma!}{(\frac{2}{\delta}+m+1)_\gamma}\right.
\nonumber \\  
& & -\left.\left(\frac{r}{R_\odot}\right)^{m\delta+2}
\,_2F_1\left(-\gamma,\frac{2}{\delta}+m;
\frac{2}{\delta}+m+1;\left(\frac{r}{R_\odot}\right)^\delta\right)
\right],
\end{eqnarray}
where $_2F_1(.)$ is Gauss' hypergeometric function and $(-
\gamma)_m=\Gamma(-\gamma+m)/\Gamma(-\gamma)$ is Pochhammer's symbol
that often appears in series expansions for hypergeometric
functions. The solution of equation (4.1) given in equation (4.2)
confirms the condition $P(r=R_\odot)=0$ and gives the central value
of the pressure according to the chosen solar model characterized
by $\delta$ and $\gamma$ in equation (2.1):
\clearpage
\begin{eqnarray}
P_c &=&\frac{9}{4\pi}G\frac{M^2_\odot}{R^4_\odot}\left[\frac{1}{\gamma!}\prod^\gamma_{i=1}\left(
\frac{3}{\delta}+i\right)\right]^2
\nonumber \\
& & \times \frac{1}{\delta}
\sum^\infty_{m=0}\frac{(-
\gamma)_m\gamma!}{m!(\frac{3}{\delta}+m)(\frac{2}{\delta}+m)
(\frac{2}{\delta}+m+1)_\gamma}.
\end{eqnarray}

\section{Temperature Distribution}
\hspace{\parindent}The simplest theory of solar structure is that
of a polytrope.
These solar models obey an equation of state of the form
$P=K\rho^{(n+1)/n}$ throughout the gas sphere. Since temperature
does not explicitly appear in this relation between $\rho$ and $P$,
equations (3.1) and (4.1) may be solved independently of the
temperature and luminosity gradients. This equation of state leads
to the Lane-Emden equation for polytropic gas spheres, which is an
ordinary differential equation of second order, but can be reduced, by
suitable transformations of the variables, to an equation of the
first order (Chandrasekhar, 1939; Horedt, 1990). In the sun
the density is so low that at the temperatures involved the solar
material behaves almost as a perfect gas, having molecular weight
$\mu$ and obeying the perfect gas law
\begin{equation}
P=\frac{kN_A}{\mu}\rho T,
\end{equation}
where k is Boltzmann's constant and $N_A$ Avogadro's number.
Substituting $\rho\propto M_\odot/R_\odot^3$ and $P\propto GM_\odot^2/R_\odot^4$ in (5.1)
reveals the dependence of temperature on the mass and radius of the sun
$T\propto \mu M_\odot/R_\odot$, where the constant of proportionality depends
on the mass distribution and the radial distance. We obtain the
detailed temperature distribution throughout the sun by using
equations (2.1) and (4.2) to rewrite the equation of state given in
(5.1), that is:
\begin{eqnarray}
T(\frac{r}{R_\odot})& = &
3\frac{\mu}{kN_A}G\frac{M_\odot}{R_\odot}\left[\frac{1}{\gamma!}\prod^\gamma_{i=1}\left(\frac{3}{\delta}+i\right)\right]\nonumber \\
& & \times \frac{1}{\delta^2}
\frac{1}{[1-
(\frac{r}{R_\odot})^\delta]^\gamma}
\sum^\infty_{m=0}\frac{(-
\gamma)_m}{m!(\frac{3}{\delta}+m)(\frac{2}{\delta}+m)}
\left[\frac{\gamma!}{(\frac{2}{\delta}+m+1)_\gamma}\right.\nonumber
\\
& & \left.-
\left(\frac{r}{R_\odot}\right)^{m\delta+2}\,_2F_1\left(-\gamma,
\frac{2}{\delta}+m;
\frac{2}{\delta}+m+1;\left(\frac{r}{R_\odot}\right)^\delta\right)
\right].
\end{eqnarray}
Equation (5.2)
satisfies the boundary condition $T(r=R_\odot)=0$ and allows to
determine the central value of temperature of the sun as a function
of $\delta$ and $\gamma$ contained in equation (2.1):
\begin{eqnarray}
T_c &=&3\frac{\mu}{kN_A}G\frac{M_\odot}{R_\odot}\left[\frac{1}{\gamma!}\prod^\gamma_{i=1}\left(
\frac{3}{\delta}\right)\right]\nonumber \\
& & \times\frac{1}{\delta^2}
\sum^\infty_{m=0}\frac{(-
\gamma)_m\gamma!}{m!(\frac{3}{\delta}+m)(\frac{2}{\delta}+m)
(\frac{2}{\delta}+m+1)_\gamma}.
\end{eqnarray}
The procedure is to construct a model for the
internal structure of the sun by assuming a matter density
distribution and subsequently integrating the system of
differential equations. At this point  two remarks are in order. The contributions
of the radiation pressure to the total pressure and the radial
dependence of the mean molecular weigth $\mu$ have been neglected
in equation (5.1).\par
The total pressure $P$ at any point in the sun is the sum of the
gas pressure and the radiation pressure, $P=P_g+P_r$, where $P_g$
is given in equation (5.1) and $P_r=\frac{1}{3}aT^4$, where $a$ is
a
constant. Writing $P_g=\beta P$ and hence $P_r=(1-\beta)P$, it
follows that $P=aT^4/3(1-\beta).$ This ratio of radiation pressure
to gas pressure increases towards the center of the Sun but even
there the gas pressure exceeds the radiation pressure by more than
three orders of magnitude. This justifies neglectingt the radiation
pressure has been neglected in equation (5.1) (Chandrasekhar,
1939).\par
The outward flow of energy inside the sun is driven by
the temperature gradient and resisted by the opacity of the
material. The nuclear energy generated within the sun has to
continually replenish that radiated away from the surface. This
energy generation by nuclear reactions causes the solar chemical
composition to change and keeps the sun evolving. Since the gas
inside the sun is completely ionised, the mean molecular
weight $\mu$ in equation (5.1) is given by
$\mu=(2X+\frac{3}{4}Y+\frac{1}{2}Z)^{-1},$ where $X, Y, Z$ are
relative abundances by mass of hydrogen, helium, and heavy
elements $(X+Y+Z=1)$. The dependence of $X, Y, Z$ on the radial
distance
variable $r$, which is governed by kinetic equations, cannot be
determined in the procedure of constructing an analytic solar model
by assuming a matter density distribution as given in equation
(2.1).
Thus, the mean molecular weigth has to be treated as constant in
the
following. This assumption does not reflect the situation in
the real sun because nuclear reactions have changed the originally
uniform chemical composition throughout the sun.

\section{Nuclear Energy Generation Rate}\par
\hspace{\parindent}Nuclear energy production in the sun depends
heavily on the
temperature of the material and is very concentrated towards the
center of the sun. This explains that
calculations of the internal structure of solar-type-stars made
considerable progress even before the physical mechanism of the
production of energy by nuclear reactions was understood. The rate
of nuclear energy generation can be written (Mathai and Haubold
(1988), as
\begin{equation}
\epsilon(\rho, T)=\epsilon_0 \rho^n(r)T^m(r),
\end{equation}
where $\epsilon_0$ is a physical constant depending only on the
chemical composition of the solar material and the units chosen.
Substituting $\rho(r)$ and $T(r)$ in equation (6.1) by
equations (2.1) and (5.2), respectively, and taking advantage of
the
fact that $0\leq(\frac{r}{R_\odot})\leq1,$ the energy generation
rate
in equation (6.1) can be represented in the form of a polynomial
\begin{equation}
\epsilon\left(\frac{r}{R_\odot}\right)=\epsilon_0\rho^n_c T^m_c f
\left(\frac{r}{R_\odot}\right)^{\delta
s+2q+\delta[n_1+2n_2+\ldots+(2\gamma)n_{2\gamma}]},
\end{equation}
where $f$ denotes the expression
$$f(\delta,
\gamma,m,n;q;a_0,a_1,\ldots,a_
{2\gamma})$$
\begin{eqnarray}
= & {\displaystyle\sum^{\gamma(m-n)}_{s=0}} & \frac{[\gamma(m-
n)]_s}{s!}\sum^m_{q=0}\frac{(-
m)_q}{q!}\left(\frac{1}{\eta(\gamma)}\right)^q\nonumber \\
& &
\times\sum^q_{n_0=0}\ldots\sum^q_{n_{2\gamma}=0}\frac{q!a_0^{n_0
}a_1^{
n_1}\ldots
a_{2\gamma}^{n_{2\gamma}}}{n_0!n_1!\ldots, n_{2{\gamma!}}},
\end{eqnarray}
$$n_0+n_1+\ldots+n_{2\gamma}=q,$$
and
\begin{equation}
\eta(\gamma)=\sum^\gamma_{\nu=0}\frac{(-
\gamma)_\nu}{\nu!}\frac{1}{(\frac{2}{\delta}+\nu)(\frac{3}
{\delta}+\nu)}\frac{\gamma!}{(\frac{2}{\delta}+\nu+1)_\gamma}.
\end{equation}
Note that the representation of $\epsilon(r/R_\odot)$ in equation
(6.2) is essentially determined by the four free parameters
$\delta,
\gamma, n$ and m in equation (2.1) and (6.1), the only
restriction being 
that $\gamma$ be a positive integer. The coefficients $a_0, a_1,
\ldots a_{2\gamma}$ in equation (6.3) are determined by the
following polynomial of degree $2\gamma$ in
$\left(\frac{r}{R_\odot}\right)^\delta:$
\begin{eqnarray}
\sum^\gamma_{m_1=0}\sum^\gamma_{m_2=0} & {\displaystyle\frac{(-
\gamma)_{m_1}}{m_1!}} & \frac{(-
\gamma)_{m_2}}{m_2!}\frac{1}{\left(\frac{2}{\delta}+m_1\right)
\left(\frac{3}{\delta}+m_1\right)\left(\frac{2}{\delta}+m_1+m_2
\right)}\nonumber \\
& & \times
\left[\left(\frac{r}{R_\odot}\right)^\delta\right]^{m_1+m
_2}
=\sum^{2\gamma}_{m_3=0}a_{m_3}\left[\left(\frac{r}{R_\odot}\right
)^\delta\right]^{m_3}.
\end{eqnarray}
\section{Luminosity Function}
\hspace{\parindent}Let $L(r)$ be the function representing the flow
of integrated
radiation across a sphere of radius $r$. If $\epsilon$ is the
energy
produced per unit time by nuclear reactions in each unit mass of
solar material, then the balance between energy generation inside the sun and energy loss through its surface is governed by
the equation of energy conservation
\begin{equation}
\frac{dL(r)}{dr}=4\pi r^2\rho(r)\epsilon (r),
\end{equation}
where $\epsilon (r)$ is given by equation (6.2). Integrating
equation (7.1) over the sun's interior leads to the luminosity
function in terms of Gauss' hypergeometric function:
\begin{eqnarray}
L\left(\frac{r}{R_\odot}\right) &=& 4\pi
\epsilon_0 \rho_c^{n+1}T_c^mR^3_\odot\nonumber \\
& &\times\frac{1}{\delta}f\frac{1}{s^*}
\left(\frac{r}{R_\odot}\right)^{\delta s^*}\,_2F_1\left(-
\gamma,s^*;s^*+1;\left(\frac{r}{R_\odot}\right)^\delta\right),
\end{eqnarray}
$$s^*=s+\frac{1}{\delta}(3+2q)+n_1+2n_2+\ldots+(2\gamma)n_
{2\gamma},$$
where $f$ is given in equation (6.3) with s substituted by $s^*.$
Equation (7.2) satisfies the condition $L(r=0)=0$ and gives for the
total energy output $L(r=R_\odot)=L_\odot,$
\begin{equation}
L_\odot=4\pi \epsilon_0 \rho_c^{n+1}T_c^mR_\odot^3 \frac{1}{\delta}
f\frac{\gamma!}{\prod^\gamma_{i=0}(s^*+i)}.
\end{equation}
\section{Conclusions}
\hspace{\parindent}Assuming a two-parameter family of matter
density distributions in
equation (2.1) made possible the analytic integration of the
differential equations for conservation of mass, momentum, and
energy throughout the sun, respectively equations (3.1), (4.1) and (7.1).
This procedure shows that hydrostatic equilibrium, equation of
state, and overall energy conservation determine the state of the
central solar conditions, particularly the gravitationally
stabilized solar fusion reactor. The mathematical method chosen
reveals the factors of proportionality which are kept undetermined
in dimensional analysis, commonly pursued to understand
astrophysical relationships between global parameters of the sun.
A common mathematical element of the derived distributions of mass
(3.2), pressure (4.2), temperature (5.2), and luminosity (7.2)
throughout the
sun is Gauss' hypergeometric function $_2F_1(\cdot)$ which is
numerically easily accessible through mathematical programs like Mathematica (Wolfram, 1993).\par
It has been emphasized above that the assumption of an analytic
matter density distribution means that the equation for the transport of
energy by radiation through the sun can be satisfied at only one
specific 
point in the sun. The flow of radiant energy per unit area through
the sun is proportional to the ratio of radiation pressure gradient
and opacity per unit volume. That is
\begin{equation}
H\propto\frac{d(\frac{1}{3}aT^4)/dr}{\kappa \rho}\propto
\frac{T^3dT/dr}{\kappa\rho},
\end{equation}
where $\kappa$ denotes the opacity per mass unit at temperature $T$
and density $\rho$. Because the energy flowing out through the sun
is transported by radiation, we find for the luminosity $L_\odot$ 
\begin{equation}
L_\odot\propto4\pi R_\odot^2 H \propto \frac{R_\odot^2 T^3 dT/dr}{\kappa
\rho}.
\end{equation}
Since for a given solar structure $\rho\propto M_\odot/R_\odot^3$
and $T\propto M_\odot/R_\odot$, it follows for $L_\odot$ that
\begin{equation}
L_\odot \propto \frac{1}{\kappa}M_\odot^3.
\end{equation}
For solar the composition, Kramer's power law approximation for
the opacity given by
\begin{equation}
\kappa \propto \kappa_0 \rho T^{-7/2},
\end{equation}
where $\kappa_0$ is a physical constant depending on the chemical
composition of the solar material and the units chosen, which leads
to a luminosity-mass-radius relation,
\begin{equation}
L_\odot \propto M_\odot^{11/2}R_\odot^{-1/2}.
\end{equation}
The differential equation governing the outward flow of energy
driven by the temperature gradient and resisted by opacity in (8.1)
is
\begin{equation}
L_\odot=4\pi r^2H=-\left(\frac{16\pi ac}{3\kappa
\rho}\right)r^2T^3\frac{dT}{dr},
\end{equation}
where $c$ denotes the velocity of light. To satisfy the radiative
energy transport, taking into account the density
distribution assumed in equation (2.1), the temperature gradient
$dT/dr$ in equation (8.6) has to be equal to the temperature
gradient
in equation (5.2). This condition can be satisfied at only one
specific point in the solar interior, for example at the boundary
of the nuclear energy producing core region (at $r\approx
0.3R_\odot$ where $L\approx L_\odot).$\par
\bigskip
\begin{center}
Acknowledgements
\end{center}
\noindent
I would like to dedicate this paper to Professor A.M. Mathai, who
understands that the study of
theoretical astrophysics is an obsession, not a profession.
Hans and Barbara Haubold would like to thank him for the long
lasting support and encouragement.\par
\clearpage
\noindent
\begin{tabbing}
[122] \=  \kill 
[1]\> Abdurashitov, J.N., et al. (1994) Results from SAGE \\
\>(The Russian-American Gallium solar neutrino Experiment),\\
\> Physics Letters, \underline{B 328}, 234-248.\\[0.1cm]
[2] \>Anselmann, P., et al. (1994) GALLEX results from the first 30 solar\\ 
\>neutrino runs, Physics Letters, \underline{B 327}, 377-385.\\[0.1cm]
[3]\> Bahcall, J.N. (1989) {\it Neutrino Astrophysics}, Cambridge
University Press,\\
\> New York.\\[0,1cm]
[4] \>Bhaskar, R. and Nigam, A. (1991) Qualitative explanations
of red giant\\ 
\>formation, Astrophys. J., \underline{372}, 592-596.\\[0,1cm]
[5] \>Chandrasekhar, S. (1939/1957) {\it An Introduction to the
Study
of}\\
\>{\it Stellar Structure}, Dover Publications, Inc., New
York.\\[0,1cm]
[6] \>Davis Jr., R. (1993) Report on the Homestake solar neutrino experiment,\\ 
\>in {\it Frontiers of Neutrino Astrophysics}, Y. Suzuki and K. Nakamura, eds.,\\ 
\>Universal Academy Press, Inc., Tokyo.\\[0.1cm]
[7] \>Horedt, G.P. (1990) Analytical and numerical values of
Emden-\\
\>Chandrasekhar associated functions, Astrophys. J.,
\underline{357}, 560-565.\\[0,1cm]
[8] \>Luke, Y.L. (1969) {\it The Special Functions and their
Approximations},\\
\> Vols. I and II, Academic Press, New
York.\\[0,1cm]
[9] \>Mathai, A.M. (1993) {\it A Handbook of Generalized Special
Functions}\\
\>{\it for Statistical and Physical Sciences}, Clarendon Press,
Oxford.\\[0,1cm]
\end{tabbing}
\clearpage
\begin{tabbing}
[122] \= \kill
[10] \>Mathai, A.M. and Haubold, H.J. (1988) {\it Modern Problems in
Nuclear and}\\ 
\>{\it Neutrino Astrophysics}, Akademie-Verlag, Berlin.\\[0,1cm]
[11]\>Nakamura, K. (1993) Recent results from Kamiokande solar neutrino observations,\\ \>Nuclear Physics, \underline{B 31}, 105-110.\\[0.1cm]
[12] \>Noels, A., Papy, R., and Remy, F. (1993) Visualizing the
evolution of a star\\ 
\> on a graphic workstation, Computers in Physics,
\underline{7}, No.2, 202-207.\\[0,1cm]
[13]\>Stein, R.F. (1966) Stellar evolution: A survey with analytic
models, in {\it Stellar}\\ 
\> {\it Evolution}, R.F. Stein and A.G.W. Cameron, eds., Plenum
Press, New York.\\[0,1cm]
[14] \>Wolfram, S. (1991) {\it Mathematica - A System for Doing
Mathematics by}\\ 
\>{\it Computer}, Addison-Wesley Publishing Company Inc., Redwood
City.\\
\end{tabbing}
\end{document}